\RequirePackage{ifpdf}
\ifpdf 
\documentclass[pdftex]{sigma}
\else
\documentclass{sigma}
\fi

\begin{document}

\allowdisplaybreaks

\renewcommand{\thefootnote}{$\star$}

\renewcommand{\PaperNumber}{006}

\FirstPageHeading

\ShortArticleName{Heisenberg-Type  Families  in $U_q( {\widehat{sl_2} })$ }

\ArticleName{Heisenberg-Type  Families  in $\boldsymbol{U_q({\widehat{sl_2}})}$\footnote{This paper is a contribution to the Proceedings of the XVIIth International Colloquium on Integrable Systems and Quantum Symmetries (June 19--22, 2008, Prague, Czech Republic). The full collection
is available at
\href{http://www.emis.de/journals/SIGMA/ISQS2008.html}{http://www.emis.de/journals/SIGMA/ISQS2008.html}}}

\Author{Alexander ZUEVSKY}

\AuthorNameForHeading{A. Zuevsky}

\Address{Max-Planck Institut f\"ur Mathematik, Vivatsgasse 7, 53111 Bonn, Germany}
\Email{\href{mailto:zuevsky@mpim-bonn.mpg.de}{zuevsky@mpim-bonn.mpg.de}}

\ArticleDates{Received October 20, 2008, in f\/inal form January 13,
2009; Published online January 15, 2009}

\Abstract{Using the second Drinfeld formulation of the quantized universal enveloping algebra
$U_q(\widehat{sl_2})$ we introduce a family of its Heisenberg-type elements which are endowed with a
 deformed commutator and satisfy properties similar to
 generators of a Heisenberg subalgebra.
Explicit expressions for new family of generators are found.}


\Classification{17B37; 20G42; 81R50}

\section{Introduction}
The purpose of this paper is to introduce a new family of elements
in the quantized enveloping algebra $U_q(\widehat{sl_2})$ of an af\/f\/ine Lie algebra $\widehat{sl_2}$.
This Heisenberg-type families possess properties similar to ordinary Heisenberg algebras.
 Heisenberg subalgebras of af\/f\/ine Lie algebras and of their $q$-deformed enveloping algebras
are being actively used in various domains of mathematical physics.
 Most important applications of Heisenberg subalgebras can be found in the f\/ield of classical
and quantum integrable models and f\/ield theories.
  Vertex operator constructions both for af\/f\/ine Lie algebras \cite{kac}
and for $q$-deformations of their universal enveloping algebras~\cite{jm, savzu} are essentially based on Heisenberg subalgebras.
  Given a quantized universal enveloping algebra~$U_q(\widehat{\cal G})$ of an af\/f\/ine Kac--Moody Lie
algebra~$\widehat{\cal G}$, it is rather important to be able to extract explicitly generators of a
Heisenberg subalgebra associated to a chosen grading of $U_q(\widehat{\cal G})$ which is not always a
trivial task.
 For instance, one can easily recognize elements of a Heisenberg subalgebra
 among generators in the homogeneous grading of the second Drinfeld realization of
$U_q(\widehat{sl_2})$ \cite{jm} while it is not obvious
 how to extract a Heisenberg subalgebra in the principal grading.
Ideally, one would expect to obtain a realization of the Heisenberg subalgebra
associated to the principal graiding of $U_q(\widehat{sl_2})$ which would involve
ordinary (rather then $q$-deformed) commutator in commutation relations
with certain elements in the family giving central elements.
 This would lead to many direct applications both in quantum groups and
quantum integrable theories in analogy with the homogeneous graiding case.
In \cite{enr} the principal commuting subalgebra in the nilponent part of $U_q(\widehat{sl_2})$
was constructed. Its elements expressed in $q$-commuting coordinates
commute with respect to the $q$-deformed bracket.

In this paper we introduce another possible version of a family elements in $U_q(\widehat{sl_2})$
which could play a role similar to ordinary Heisenberg subalgebra.
Our general idea is to form certain sets of $U_q(\widehat{sl_2})$-elements containing linear
combinations of generators $x^{\pm}_n$, $n \in \mathbb{Z}$, multiplied by various powers of $K$
and the central element $\gamma$. Under certain conditions on corresponding powers
we obtain commutation relations for a Heisenberg-type family with respect to
 an integral $p$-th power of $K \in U_q(\widehat{sl_2})$-deformed commutator.
We consider this as some further generalization of various $q$-deformed commutator
algebras (in particular, $q$-bracket Heisenberg subalgebras)
 which f\/ind numerous examples in quantum algebras and applications in integrable models.
Though the commutation relation we use in order to def\/ine a Heisenberg-type family
look quite non-standard we believe that these families and their properties deserve
a consideration as a new structure inside the quantized universal enveloping algebra
of $\widehat{sl_2}$ even when $p \ne 0$.

The paper is organized as follows. In Section~\ref{section2} we recall the def\/inition of the quantized universal
enveloping algebra  $U_q(\widehat{sl_2})$ in the second Drinfeld realization.
In Section~\ref{section3} we f\/ind explicit expressions for elements of Heisenberg-type families.
Then we prove their commutation relations.
We conclude by making comments on possible generalizations and applications.

\section[Second Drinfeld realization of $U_q(\widehat{sl_2})$]{Second Drinfeld realization of $\boldsymbol{U_q(\widehat{sl_2})}$}\label{section2}

Let us recall the second Drinfeld realization  \cite{drinf1, frenk} of the quantized universal
enveloping algebra $U_q(\widehat{sl_2})$. It is generated by the elements $\{ x^{\pm}_k, \, k \in
\mathbb{Z}; \, a_n, \,  n \in \{ \mathbb{Z} \backslash 0 \}; \,  \gamma^{ \pm \frac{1}{2} },\,
K \}$, subject to the commutation relations
\begin{gather}
[K,a_k]=0, \qquad K   x^{\pm}_k   K^{-1} =q^{\pm 2}   x^{\pm}_k,
\qquad
[a_k, a_l] =\delta_{k,-l}   \frac{[2k]}{k}
\frac{\gamma^k -\gamma^{-k} }{q-q^{-1}},
\nonumber\\
[a_n, x^{\pm}_k ]= \pm  \frac{[2n]}{n}   \gamma^{\mp \frac{|n|}{2}}   x^{\pm}_{n+k},\label{2d}
\qquad
[x^+_n, x^-_k]=  \frac{1} {q-q^{-1}}
 \big(
\gamma^{\frac{1}{2}(n-k)}  \psi_{n+k}
-\gamma^{-\frac{1}{2}(n-k)}  \phi_{n+k}
\big),
\\
x^\pm_{k+1}   x^\pm_l - q^{\pm 2}  x^\pm_l   x^\pm_{k+1}=
q^{\pm 2}  x^\pm_k   x^\pm_{l+1} - x^\pm_{l+1}   x^\pm_k,\nonumber
\end{gather}
 where $\gamma^{\pm \frac{1}{2} }$ belong to the center of
$U_q(\widehat{sl_2})$, and
\[
[n]\equiv \frac{q^n - q^{-n}}{q-q^{-1}}.
\]
The elements $\phi_k$ and $\psi_{-k}$, for non-negative integers $k \in \mathbb{Z}_+$,  are related
to $a_{\pm k}$ by means of the expressions
\begin{gather}
\sum\limits_{m=0}^{\infty}
\psi_m   z^{-m}= K\exp\left( (q-q^{-1} ) \sum\limits_{k=1}^{+\infty} a_{k}   z^{-k} \right),
\nonumber\\
\sum\limits_{m=0}^{\infty}
\phi_{-m}   z^{m}= K^{-1}\exp\left( -(q-q^{-1} )
\sum\limits_{k=1}^{+\infty} a_{-k}    z^{k}
 \right),\label{psi}
\end{gather}
i.e.,
$\psi_m=0$,   $m < 0,$
$\phi_m=0,$
$m > 0$.

\section{Heisenberg-type families}\label{section3}
In this section we introduce a family of $U_q(\widehat{sl_2})$-elements
which have properties similar to a~ordinary Heisenberg subalgebra of
an af\/f\/ine Kac--Moody Lie algebra \cite{kac}.
We consider families of linear combinations of $x^{\pm}_n$-generators of $U_q(\widehat{sl_2})$
 multiplied by powers of $K$ and central element~$\gamma$.

Let us introduce for $m$, $l$, $\eta$, $\theta \in \mathbb{Z}$, $n \in \mathbb{Z}_+$,
the following elements:
\begin{gather}
\label{e2dp}
E_n^\pm(m, \eta)= \gamma^{\pm(n+ \frac{1}{2})}    x^+_n   K^m +  x^-_{n+1}    K^\eta,
\\
\label{e2dm}
E_{-n-1}^\pm (l, \theta)=    x^+_{-n-1}   K^l +  \gamma^{\pm(n+ \frac{1}{2})} x^-_{-n}    K^\theta.
\end{gather}
Denote also for some $p\in \mathbb{Z}$, a deformed commutator
\begin{equation}
\label{defcomm}
[A,B]_{K^p}=A  K^p  B - B K^p  A.
\end{equation}

We then formulate

\medskip

\noindent
{\bf Proposition.} {\it Let $p$, $m \in \mathbb{Z}$, $l=m$, $\theta=\eta =-m-2p$.
Then the family of elements
\begin{equation}
\label{family}
 \{
  E_{p, n}^{\pm}(m),
 \,
 E_{p, -n-1}^{\pm} (m),
\, n \in \mathbb{Z}_+ \},
\end{equation}
where
$E_{p, n}^{\pm} (m) \equiv E_{n}^\pm(m,  -m - 2p)$,
$E_{p, -n-1}^{\pm} (m) \equiv E_{-n-1}^\pm (m, -m - 2p)$,
we have denoted in~\eqref{e2dp},~\eqref{e2dm},
are subject to the commutation relations with $k \in \mathbb{Z}_+$,
\begin{gather}
\label{ep}
\big[ E_{p, n}^+(m), E_{p, -k-1}^+(m) \big]_{K^p}=0,\qquad {\rm for \  all} \ n<k,
\\
\label{em}
\big[ E_{p, n}^-(m) , E_{p, -k-1}^-(m) \big]_{K^p}=0,\qquad {\rm for \  all} \  n>k,
\\
\label{commc}
\big[ E_{\pm 1, n}^{ \pm}(m) , E_{\pm 1, -n-1}^{\pm}(m)  \big]_{K^{\mp 1}}=
   c_n^\pm(m),
\end{gather}
where
\begin{gather*}
 c_{n}^+(m)=  \frac{q^{-2(m-1)}}{q-q^{-1}}  \gamma^{2n+1} (\gamma^{n}-\gamma^{-n-1}),
\qquad
 c_{n}^-(m)=  \frac{q^{-2(m+1)}}{q-q^{-1}}   \gamma^{-n-1} (\gamma^{2n+2}-\gamma^{-n}),
\end{gather*}
belong to the center ${\cal Z}(U_q(\widehat{sl_2}))$ of $U_q(\widehat{sl_2})$.
}

\medskip

 We call a subset (\ref{family}) of $U_q(\widehat{sl_2})$-elements
with all appropriate $m$, $p \in \mathbb{Z}$, $n \in \mathbb{Z}_+$, such that it
satisf\/ies the commutation relations (\ref{ep})--(\ref{commc}) the {\it Heisenberg-type} family.
In particular, when $p=0$, (\ref{ep})--(\ref{em}) reduce to ordinary commutativity conditions.
Note that if we formally substitute $n \mapsto -n -1$, then,
\[
E_{p, n}^\pm (m, \eta) \mapsto \gamma^{\mp(n + 1/2)} E_{p, -n-1}^\mp(m, \eta).
\]
Under the action of an automorphism $\omega$ of $U_q(\widehat{sl_2})$ which maps
$K \mapsto K^{-1}$, $\gamma \mapsto \gamma^{-1}$,
  $x_{n}^{\pm} \mapsto x_{-n}^{\mp}$, $a_{n} \mapsto a_{-n}$, one has
$ \omega(E_{p, n}^\pm (m)) =  E_{p, -n-1}^\mp(m) \; K^{2p}$.

\noindent
\begin{proof}  The proof is the direct calculation of the commutation relations
(\ref{ep})--(\ref{commc}).
 Indeed, consider $K^p$-deformed commutator (\ref{defcomm}) of $E^{\pm}_n(m,\eta)$ and
 $E^{\pm}_{-k-1}(l, \theta)$ with some
 $m, \eta, l, \theta \in \mathbb{Z}$, $n, k \in \mathbb{Z}_+$.
Using the commutation relations (\ref{2d})
we obtain
\begin{gather*}
\left[  E_n^{\pm}(m, \eta),  E_{-k-1}^{\pm}(l, \theta) \right]_{K^p}
=\gamma^{\pm(n+k+1)}   \big(q^{ -2m-2p}
  x_n^+   x_{-k}^-  - q^{2\theta + 2p  }   x_{-k}^-   x_n^+ \big)
  K^{m+\theta + p}
\\
 \qquad {} + \big( q^{ 2\eta + 2p }    x_{n+1}^-   x_{-k-1}^+ - q^{-2l - 2p}    x_{-k-1}^+   x_{n+1}^-
\big)   K^{\eta + l+p}
\\
\qquad {} + \gamma^{\pm(n+1/2)} \big(q^{ 2m+2p}   x_n^+   x_{-k-1}^+ - q^{2l+2p}  x_{-k-1}^+   x_n^+  \big)
   K^{\eta + \theta + p}
\\
\qquad {}+ \gamma^{\pm(k+1/2)} \big(q^{ -2\eta -2p}  x_{n+1}^-   x_{-k}^-
 - q^{-2\theta -2p }   x_{-k}^-   x_{n+1}^-  \big)
    K^{\eta + \theta + p}.
\end{gather*}
Then for $m=l$, $\theta=\eta=-m-2p$, from (\ref{2d}) it follows
\begin{gather*}
\big[ E_{p, n}^{\pm} (m),  E_{p, -k-1}^{\pm} (m)\big]_{K^p}
= \frac{q^{-2(m+p)}}{q - q^{-1}}   \big[  \gamma^{\pm(n+k+1)}
\big( \gamma^{1/2(n+k)} \psi_{n-k}  - \gamma^{-1/2(n+k)} \phi_{n-k} \big)
\nonumber\\
\qquad {}- \big( \gamma^{-1/2(n+k+2)} \psi_{n-k}  - \gamma^{1/2(n+k+2)} \phi_{n-k} \big) \big]
K^{p}.
\end{gather*}
Since for $n<k$, $\psi_{n-k}=0$, and the reaming terms containing $\phi_{n-k}$ cancels,
 we obtain (\ref{ep}).
 Similarly, for $n>k$, $\phi_{n-k}=0$, the terms containing $\psi_{n-k}$ cancels,
and (\ref{em}) follows.

From (\ref{psi}) we see that $\psi_0=K$, $\phi_0=K^{-1}$. Taking $K^{\mp 1}$-deformed
commutators (\ref{defcomm}) of $E_{1, n}^{\pm}(m)$, $E_{1,-k-1}^{\pm}(m)$, we then have
\begin{gather*}
\big[ E_{1,n}^{+}(m),  E_{1, -n-1}^{+}(m)  \big]_{K^{-1}}
= \frac{q^{-2m+2}}{q - q^{-1}}  \big[ \gamma^{2n+1}  (\gamma^{n} K )\;K^{-1}
-
  (\gamma^{-n-1} K ) K^{-1}\big]
=c_{n}^+(m),
\\
\big[ E_{-1, n}^{-}(m),  E_{-1, -n-1}^{-}(m)  \big]_K
= \frac{q^{-2m-2}}{q - q^{-1}} \big[\gamma^{-2n-1} K (  - \gamma^{-n}K^{-1})
+
K( \gamma^{n+1}K^{-1}) \big]
=c_{n}^-(m),
\end{gather*}
for $n=k$.
\end{proof}

\section{Conclusions}
In the second Drinfeld realization of $U_q(\widehat{sl_2})$ we have def\/ined a subset of elements
that constitutes a Heisenberg-type family, explicitly constructed their elements, and proved
corresponding commutation relations.
 Properties of a Heisenberg-type family are similar to ordinary Heisenberg subalgebra properties.
These families might be very useful in construction of special types of vertex operators
in $U_q(\widehat{sl_2})$,
 and, in particular,
might have their further applications in the soliton theory
of non-linear integrable partial dif\/ferential equations \cite{savzu}.
 One of our aims to introduce Heisenberg-type families is the
development of corresponding vertex operator representation which plays the main
role in the theory of quantum soliton operators in exactly solvable f\/ield models associated to
the inf\/inite-dimensional Lie algebra $\widehat{sl_2}$ \cite{savzu}.

Finally, we would like also to make some comments comparing present work to \cite{enr}
where the quantum principal commutative subalgebra in $U_q(\widehat{sl_2})$
associated to the principal grading of~$\widehat{sl_2}$ was found.
Here we introduce Heisenberg-type families of $U_q(\widehat{sl_2})$
in the principal grading of~$U_q(\widehat{sl_2})$~\cite{savzu}.
Although we use $K^p$-deformed commutators (which for $p=1$ can be seen quite similar to $q$-deformed
commutators in~\cite{enr}) these two approaches are quite dif\/ferent.
We prefer to work with the explicit set of~$U_q(\widehat{sl_2})$ generators
($q$-commutative coordinates) in its second Drinfeld realization \cite{drinf1},
 and introduce elements of our Heisenberg-type families not involving lattice constructions or
 trace invariants as in~\cite{enr}.
 A generalization of our results to an arbitrary~$\widehat{\cal G}$ case does not
 face any serious technical problems.
We assume that Heisenberg-type families introduced which exhibit properties similar
to a Heisenberg subalgebra in $U_q(\widehat{sl_2})$ are not the most general ones.
 At the same time the construction described in this paper allows further generalization to
cases of arbitrary af\/f\/ine Lie algebras \cite{sln}.
Using formulae from \cite{enr} we see that even in the $q$-commutator case
there exist more complicated $q$-commutative elements in $U_q(\widehat{sl_2})$. Thus one
would expect the same phenomena for $K^p$-deformed algebras.

More advanced examples of Heisenberg-type families associated to various gradings
of $U_q(\widehat{\cal G})$ in the Drinfeld--Jimbo and second Drinfeld realizations as well as
corresponding vertex operators will be discussed in a forthcoming paper \cite{sln}.

\subsection*{Acknowledgements}
We would like to thank A.~Perelomov, D.~Talalaev and M.~Tuite for illuminating discussions and comments.
  Making use of the occasion, the author would like to express his gratitude to
 the Max-Planck-Institut f\"ur Mathematik in Bonn where this work has been completed.

\pdfbookmark[1]{References}{ref}
\LastPageEnding
\end{document}